# On traces of $d$-stresses in the skeletons of lower dimensions of homology $d$-manifolds

R. M. Erdahl[*]
erdalhr@post.queensu.ca
Department of Mathematics and Statistics, Queen's University
K. A. Rybnikov,[†]
rybnikov@mast.queensu.ca
Department of Mathematics and Statistics, Queen's University and
the Fields Institute for Research in Mathematical Sciences
S. S. Ryshkov,[‡]
sergei.s@ryshkov.mian.su
Steklov Mathematical Institute and
Moscow State University

February 1999

**Abstract**

We show how a $d$-stress on a piecewise-linear realization of an oriented (non-simplicial, in general) $d$-manifold in $\mathbb{R}^d$ naturally induces stresses of lower dimensions on this manifold, and discuss implications of this construction to the analysis of self-stresses in spatial frameworks. The constructed mappings are not linear, but polynomial. In 1860-70s J. C. Maxwell described an interesting relationship between self-stresses in planar frameworks and vertical projections of polyhedral 2-surfaces. We offer a partial analog of Maxwell correspondence for self-stresses in spatial frameworks and vertical projections of 3-dimensional surfaces based on our construction of polynomial mappings. Applying this theorem we derive a class of three-dimensional spider webs similar to the family of two-dimensional spider webs described by Maxwell. In addition, we conjecture an important property of our mappings which is supported by a heuristic count based on the lower bound

[*]the work of R. M. Erdahl was supported in part by grants from NSERC
[†]the work of K. A. Rybnikov was supported in part by Fields Graduate Scholarships
[‡]the work of S. S. Ryshkov was supported in part by grants from Russian Foundation for Fundamental Research

theorem ($g_2(d+1) = dim\ Stress_2 \geq 0$) for $d$-pseudomanifolds generically realized in $\mathbb{R}^{d+1}$ [12].

# 1 Introduction

If $(E,V)$ are the edges and the vertices of a framework (possibly infinite) in $\mathbb{R}^d$, then a self-stress (or simply stress) is an assignment of real numbers $s_{ij} = s_{ji}$ to the edges, a tension if the sign is positive or a compression if the sign is negative, so that the equilibrium conditions

$$\sum_{\{j\ |\ (ij)\in E\}} s_{ij} \frac{(\mathbf{v}_j - \mathbf{v}_i)}{\|\mathbf{v}_j - \mathbf{v}_i\|} = 0$$

hold at each vertex $\mathbf{v}_i \in \mathbb{R}^d$. In fact, the space of self-stresses is the left null-space of the *rigidity matrix* $RM(E,V)$ of the framework $(E,V)$, which is constructed as follows. Let $M$ be the incidence matrix for the edges and vertices of $(E,V)$, where the rows are indexed by the edges and the columns by the vertices. Thus $M(i,j) = 1$ if and only if $v_j \subset \partial e_i$, but is equal to 0 otherwise. The matrix $RM$ is obtained by replacing unit entries of $M$ by the corresponding inward unit normals to edges at their vertices, and zero entries by the zero vector in $\mathbb{R}^d$; these replacement vectors are taken to be row vectors. The left null-space of $RM$ (which consists of vectors having $|E|$ coordinates) is the space of self-stresses. The dimension of this space is equal to $|E|$ minus the number of independent rows of $RM$, or in other words, the row corank of $RM$. The row rank of $RM$ is exactly the dimension of the subspace of external loads that can be resolved by the framework. If all external loads can be resolved, the framework is called statically rigid. Since the dimension of the space of external loads in $\mathbb{R}^d$ is $dV - \binom{d+1}{2}$, the static rigidity implies that the dimension of the space of stresses is $E - dV + \binom{d+1}{2}$.

A *spider web* is a framework (with vertices at infinity usually allowed) which supports a self-stress which is strictly positive on all edges. Spider webs in $\mathbb{R}^2$ naturally appear from projections of convex surfaces. Planar and spatial spider webs serve as a tool for investigating various problems about dense packing of equal balls in $\mathbb{R}^2$ and $\mathbb{R}^3$ [1, 8, 9]. There are interesting applications of the theory of stresses in frameworks to physics, chemistry, and engineering (see [1, 3, 8, 9, 30]).

The notion of self-stress on a framework can be naturally generalized to $k$-stresses on cell-complexes. This generalization proves to be useful in the combinatorics and geometry of P.L.-manifolds, the rigidity theory, and the theory of Dirichlet-Voronoi diagrams. Such generalizations were offered by Lee [13], Tay et al. [23], Crapo and Whiteley [10], andRybnikov [18].

J.C.Maxwell discovered that the projection of a piecewise-linear sphere onto a plane induces a self-stress on the 1-skeleton of the projection [14, 15]. He also noticed that this relation is in some sense one-to-one. Later, Crapo and Whiteley [8, 26] proved that for piecewise-linear spheres realized in $\mathbb{R}^2$ (like on Figure 1) there is a natural linear isomorphism between the space of self-stresses on the 1-skeleton and the space of liftings (an operation inverse to the vertical projection of a spatial polyhedron onto the plane) considered up to the choice of a supporting plane for a facet of the lifting.

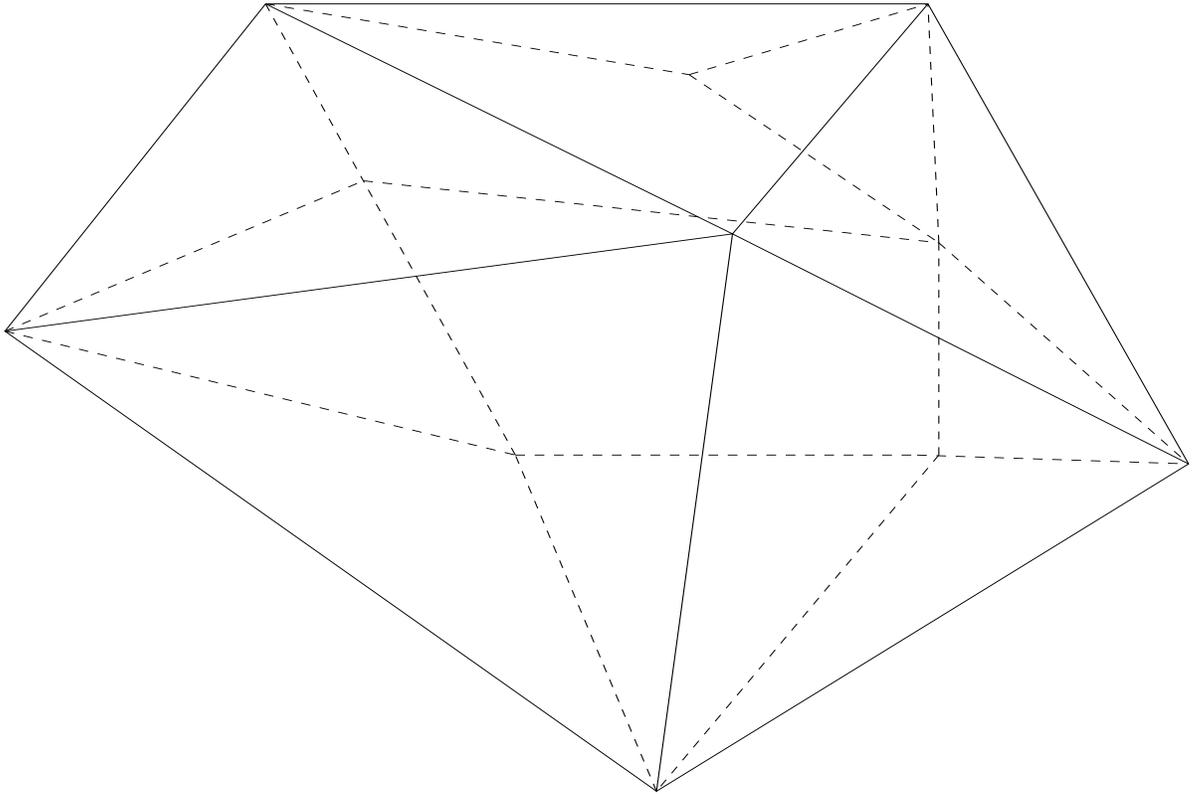

Figure 1: 2-sphere realized in $\mathbb{R}^2$

In his studies of the relationships between stresses and projections Maxwell used a new geometrical tool, so called *reciprocals*. Roughly speaking, a reciprocal for a planar piecewise-linear realization of a sphere is a special planar realization of the dual combinatorial graph of the spherical complex, such that its edges are perpendicular to the corresponding edges of the spherical complex (see Figure 2). For spheres the space of reciprocals is also isomorphic to the space of self-stresses.

As it was shown by Crapo, Whiteley and Rybnikov for certain classes of $d$-manifolds including homology spheres there is a similar connection between piecewise-linear $d$-manifolds realized in $\mathbb{R}^{d+1}$ and stresses defined on the $(d-1)$-cells of the realizations defined by the projections. In this case the equilibrium of forces is required not at the vertices, but at each $(d-2)$-cell [8, 26, 18]. Such stresses are called $d$-stresses because $d$ is the lowest dimension of a manifold for which the space of $d$-stress is non-trivial in a sense that it depends on the combinatorics of the manifold. (The space of $d$-stresses of a closed $(d-1)$-manifold realized in $\mathbb{R}^d$ is $\mathbb{R}$ or 0 depending on whether this manifold is orientable or not.)

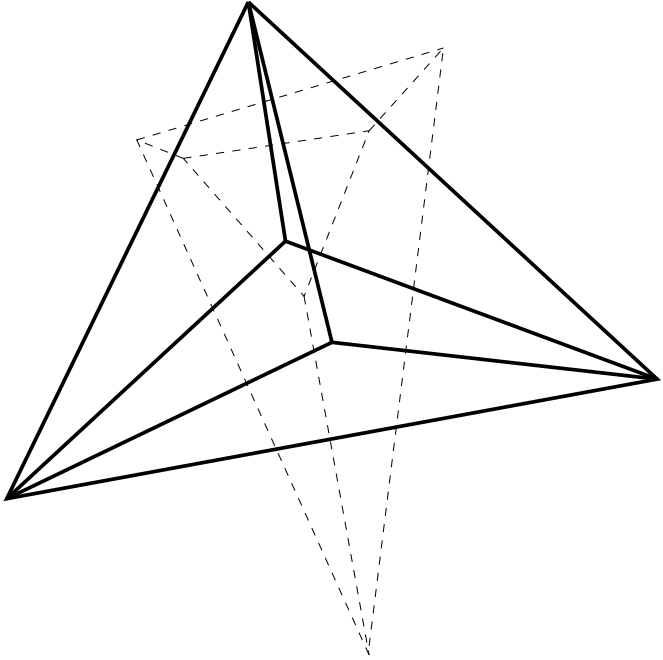

Figure 2: reciprocal for a 2-sphere in $\mathbb{R}^2$

By an informal conjecture of J. Baracs and W. Whiteley there is an analogous correspondence between projections of the boundaries of 4-polyhedra from $\mathbb{R}^4$ onto $\mathbb{R}^3$ and self-stresses in spatial frameworks [29]. This idea is motivated by Minkowski theorem on the vanishing of the sum of normals to a convex polytope at its facets (see Section 5). Since the projections of $d$-manifolds with trivial $H_1(\cdot, \mathbb{Z}_2)$ from $\mathbb{R}^{d+1}$ onto $\mathbb{R}^d$ correspond to $d$-stresses (see Section 2 and [18, 27]), one can reformulate their conjecture as the existence of a natural correspondence between $d$-stresses on a $d$-manifold realized in $\mathbb{R}^d$ and self-stresses on its 1-skeleton (in the general theory of stresses such stresses are called 2-stresses). A theorem that we prove in Section 6 can be considered as a natural one-way connection between $d$-stresses and $k$-stresses, $k < d$ on oriented piecewise-linear manifolds realized in $\mathbb{R}^d$. Therefore, in some sense our theorem supports Baracs-Whiteley hypothesis. There are other canonical mappings between the spaces of stresses of different dimensions. According to Stanley [22] and Lee [13] for Cohen-Macaulay homology $d$-spheres in $\mathbb{R}^d$ the space of $k$-stresses has the same dimension as the space of $(d-k+1)$-stresses, $k \leq \lfloor \frac{d+1}{2} \rfloor$. These isomorphisms are a very important part of the Stanley's and McMullen's proofs of the $g$-theorem [17, 22]. However these isomorphisms are linear (see [22, 13]), whereas our mappings are algebraic. In general, our mappings are *not bijective*, since for a generic realization of a simplicial sphere in $\mathbb{R}^3$ the dimension of the space of 2-stresses may exceed the dimension of the space of 3-stresses. For example, using the results of Lee [13] one can show that for a generic realization in $\mathbb{R}^3$ of the boundary of the 4-dimensional cross-polytope $dim(Stress_2) = 6$, but $dim(Stress_3) = 4$.

Our mappings are well-defined not only for simplicial manifolds, but also for cell-partitions of manifolds. We believe that these mappings are *injective* for any generic realization of an orientable simplicial $d$-manifold in $\mathbb{R}^d$. However, we are primarily interested in applications of the general theory of stresses and liftings to three dimensional manifolds and spatial spider webs. From this point of view our construction can be regarded as a natural extension of Maxwell correspondence between stresses and liftings to spatial frameworks. In the last paragraph of this section we sketch the main ideas and concepts employed in our construction.

Let $\Delta$ be a homology $d$-manifold decomposed into cells, each of which being a simplicial star (see [16, 21] for discussion of such pseudo-dissections). The construction of our mapping can be divided into several steps. On the first step we establish a natural one-to-one correspondence between $d$-stresses and *reciprocals* (see Section 4). A concept of reciprocal basically generalizes the notion of Maxwell reciprocal (see above). In particular, a segment whose ends are the vertices of the reciprocal corresponding to two adjacent $d$-cells of $\Delta$ is perpendicular to their common facet. This one-to-one correspondence can be established only under certain homological restrictions on the manifold, for example $H_1(\Delta, \mathbb{Z}_2) = 0$. Notice that the star of a cell in a manifold satisfies this condition. Given a $d$-stress $s$, we can construct the corresponding reciprocal $R(s, v)$ for the star of each vertex $v$ of $\Delta$. If two cells $C_1$ and $C_2$ share a face $F$, the sub-reciprocals of $R(s, C_1)$ and $R(s, C_2)$ corresponding to $F$ are congruent. Nevertheless, when $H_1(\Delta, \mathbb{Z}_2) \neq 0$, it is not possible, in general, to construct a global reciprocal (see [18]). One can consider for $\Delta$ the dual cell-decomposition $\Delta^*$. The idea of such decomposition goes back to Poincare ( for details see [16, 21]). We construct a piecewise-linear realization of a baricentric triangulation $TD$ of the dual cell-decomposition $\Delta^*$. Note that the baricentric triangulation of the original cell-decomposition of $\Delta$ is isomorphic to the baricentric triangulation of $\Delta^*$ [21]. After such dissection a cell of $\Delta^*$ can be realized in $\mathbb{R}^d$ as a simplicial star. We consider only special realizations, namely, such that the baricentric triangulation of each $k$-cell is realized in a $k$-plane. For instance, a sub-reciprocal $R(s, C)$ for the star of cell $C \subset \Delta$ defines such "flat" realization of $St(C^*)$ up to the choice of baricenters for the $k$-cells, $dim(\Delta) - dim(C) = dim(C^*) \geq k > 0$. In this case one can introduce a natural method of summation of the volumes of the oriented simplexes of the simplicial star $St(C^*)$ ($C^* \subset \Delta^*$), such that the result of the summation does not depend on the positions of baricenters of all cell of dimensions greater than 0. When $St(C^*)$ is embedded into $\mathbb{R}^d$ the result is the oriented volume of the simplicial star. That is why we call this function on "flat" realizations of oriented cells of the dual decomposition the *signed generalized volume*. Evidently, it can be equally thought of as a function on reciprocals. By a well-known Minkowski theorem the sum of outer (or inner) facet normals scaled with the $(d-1)$-volumes of facets of a $d$-simplex equals to zero. Using the orientability of $\Delta$, we will show in Section 6 that the generalized $k$-volumes of $k$-cells of $\Delta^*$ can be interpreted as the coefficients of $(d-k+1)$-stresses on $\Delta$. As it can be seen from this informal description, the main ingredients of our construction are volumes, reciprocals, duality in homology manifolds, and the notion of orientability. In fact, we suspect that our construction can be generalized for any dimension $n$, $\lfloor \frac{d+1}{2} \rfloor \leq n \leq d$ thereby providing canonic algebraic mappings from the space of $n$-stresses to the spaces of $k$-stresses, $d-n < k < n$.

## 1.1 Notation

All complexes that we consider are polyhedra (simplicial complexes) from the topological point of view. However, all theorems in this paper are stated for fixed decompositions of simplicial complexes into polyhedral *cells* (also called blocks or simplicial stars in combinatorial topology [21, 16]) which are not necessarily simplexes. We assume that all complexes have at most countable number of cells. Cells of co-dimension 1 are referred to as *facets*. We denote the star of a cell $C$ by $St(C)$, and the $k$-dimensional skeleton of a complex $\mathcal{K}$ by $Sk^k(\mathcal{K})$.

We shall consider a somewhat more general construction than an embedding or an immersion of a cell-complex into Euclidean space, such as a *piecewise-linear (PL-) realization* of a cell-complex in Euclidean space. In all geometric discussions cell-complexes will be considered as fixed piecewise-linear realizations, rather than abstract combinatorial objects. Such general construction can be helpful, for example, for studying frameworks with bar intersections, polyhedral scenes, splines over triangulations (in the planar case this point of view was adopted in [8, 23, 26]; in the three-dimensional case such PL-realizations were considered by Crapo and Whiteley in [8, 25]), and in the case of general dimension by Tay, White, and Whiteley [23]. For example, a Schlegel $d$-diagram is a PL-realization of a $(d+1)$-polytope $P$ in $\mathbb{R}^d$ obtained by radial projection of $P$ onto one of its facets.

Recall that one can identify an abstract combinatorial cell-complex $\mathcal{K}^d$ with its embedding into $\mathbb{R}^{2d+1}$ (since it can be triangulated). More formally, a PL-realization of a combinatorial simplicial complex $\mathcal{K}^d \subset \mathbb{R}^{2d+1}$ with a fixed decomposition into polyhedral cells is a continuous PL-mapping $r$ of $\mathcal{K}^d$ in $\mathbb{R}^N$ ($N \geq d$) such that *the closure of each $k$-cell, $k = 0, \ldots, d$ is embedded by $r$ into $\mathbb{R}^N$ as a "flat" (lying in a $k$-subspace) $k$-polyhedron*.

If $\Delta$ is a piecewise-linear realization of a polyhedron with a specific cell-decomposition, we shall frequently abuse notation and make no distinction between the polyhedron, its cell-decomposition and the piecewise-linear realization. If we refer to the metric, projective, or affine properties of a cell-complex, these should be understood as the properties of its fixed PL-realization. However, when we consider the combinatorial or homological properties of a cell-complex, we are referring to its abstract combinatorial structure.

A homology $k$-sphere ($k$-disk) is a polyhedron with the homology groups of a standard $k$-sphere ($k$-disk). A *homology $d$-manifold* (with boundary) is a cell-complex such that the link (in the case of a non-simplicial cell-decomposition, the link of a cell can is defined through the baricentric triangulation) of each $k$-cell, is either a homology $(d - k - 1)$-sphere or a homology $(d - k - 1)$-disk. A manifold is closed if each facet is adjacent to exactly two $d$-cells. All statements in the paper are formulated for both closed manifolds and for manifolds with a boundary, unless stated otherwise. Since we consider manifolds only from the combinatorial point of view, a manifold is always understood to be a *homology* manifold. Throughout the paper we include "good" decompositions of $\mathbb{R}^n$ (like, for example, weighted Voronoi diagrams) into the class of homology manifolds.

## 2 Stresses

The notion of a self-stress on a framework can be naturally generalized to a $k$-stress on a cell-complex of dimension at least $k - 1$. This generalization appears to be useful in the combinatorics and geometry of homology spheres, in rigidity theory, and in Voronoi's theory of parallelohedra (see [13, 18]).

Consider a piecewise-linear realization $K$ in $\mathbb{R}^N$ of a $d$-dimensional cell-complex $\mathcal{K}$. Denote by $\mathbf{n}(F, C)$ the inner unit normal to a cell $C$ at its facet $F$. $C$ need not be convex, but it is important that its boundary is an orientable closed manifold, namely sphere.

**Definition 2.1** *A real-valued function $s(\cdot)$ on the $(k - 1)$-cells of $K$ is a $k$-stress if at each internal $(k - 2)$-cell $F$ of $K$*

$$\sum_{\{C \mid F \subset C\}} s(C) vol_{k-1}(C)\, \mathbf{n}(F, C) = \mathbf{0},$$

*where the sum is taken over all $(k - 1)$-cells in the star of $F$. The quantities $s(C)$ are the coefficients of the $k$-stresses, a tension if the sign is strictly positive and a compression if the sign is strictly negative.*

It is easy to see that $k$-stresses form a linear space, and that $k$-tensions and $k$-compressions form congruent positive cones in this linear space. We denote the space of all $k$-stresses on $K$ by $Stress_k(M^d)$, the cone of all $k$-tensions by $Tension_k(M^d)$. If the coefficients $s(C)$ are not all zero the $k$-stress $s$ is called non-trivial. Figure 3 illustrates the geometry of the equilibrium condition for the 2-stress at an edge of a cell-complex in $\mathbb{R}^3$.

In the case of a stress on a framework, $s(e)$ is force per unit length, and the static force applied at the end points of edge $e$ is $s(e)\|e\|$. For a $(k-1)$-cell $C$ a $k$-stress $\mu$ is force per unit relative $(k-1)$-volume (area) of $C$, and the static force applied at a $(k-2)$-face of $C$ is $\mu\, vol_{k-1}C$. The notion of stress is well defined for fans and cell-decompositions of $\mathbb{R}^d$ with non-compact cells. In this case the volumes of $(k-1)$-cells should be leaved out of the formula, and a coefficient of stress does not have meaning of force per unit relative volume (area).

The definition of a $k$-stress can be adjusted so that the equilibrium of forces is not required at the $(k-2)$-cells adjacent to just one $(k-1)$-cell. For instance, it makes sense when one is describing connections between splines and stresses [27]. Another example is analysis of stresses in frameworks with fixed vertices. In this case the equilibrium of forces is not required at the fixed vertices (called pinned vertices in the planar case).

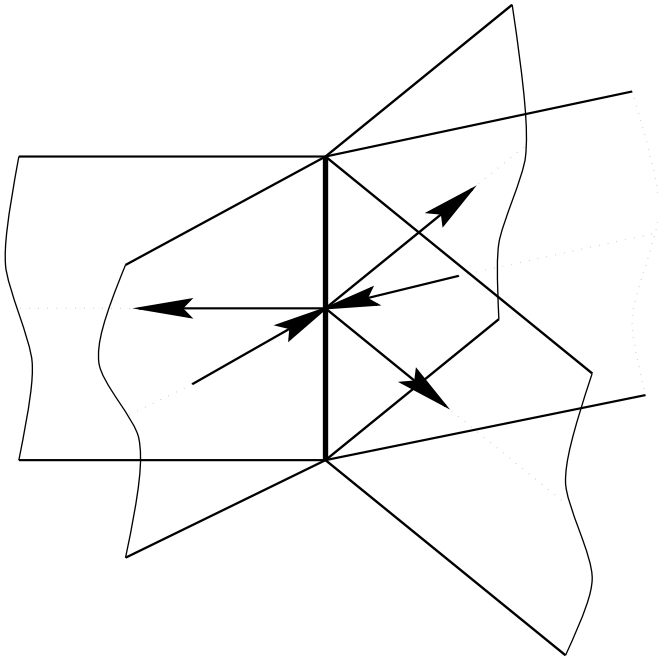

Figure 3: Equilibrium of forces at 1-cell

As in the case of frameworks, the linear space of $k$-stresses can be characterized as the left null space of a geometric matrix $RM_k$ which is constructed as follows. Let $M_k$ be the incidence matrix for the $k$- and $(k-1)$-cells of $K$, where the rows are indexed by the $k$-cells and the columns by $(k-1)$-cells. Thus $M_k(i,j) = 1$ if and only if $C_j^{k-1} \subset \partial C_i^k$, but is equal to 0 otherwise. The matrix $RM_k$ is obtained by replacing unit entries of $M_k$ by the corresponding positively oriented unit normal vectors, and zero entries by the zero vector; these replacement vectors are taken to be row vectors. The left null-space of $RM_k$ which consists of vectors (with the number of components equal to the number of $k$-cells) is the space of $(k+1)$-stresses.

The notion of $k$-stress on simplicial complexes was introduced by Lee [13]. For a simplicial complex a $k$-stress can be interpreted as an element of a certain quotient of the face-ring of the complex $K$. Let $K$ is a simplicial complex in $\mathbb{R}^d$, with vertex set $\mathbf{v}_1, \ldots, \mathbf{v}_n$. Then, in Lee's terminology the space of *affine* $k$-stresses on $K$ is the linear subspace of polynomials of degree $k$ of $R/V$, where $R$ is the Stanley-Reisner ring of $K$, and V is the ideal generated by linear forms $\sum_{i=1}^n v_{ki} x_i$ ($k = 1, \ldots, d$), and $\sum_{i=1}^n x_i$ (see [13, 23]). For a simplicial complex $K$ in $\mathbb{R}^d$ our $k$-stress on $K$ is the same as Lee's affine $k$-stress on $K$. In fact, Lee considered two types of stress: linear and affine. Lee formulated most of his theorems in terms of so-called *linear* stresses. For generic realizations of $K$ in $\mathbb{R}^d$ the space of our $k$-stresses is isomorphic to the space of Lee's *linear* $k$-stresses for $K$ realized generically in $\mathbb{R}^{d+1}$. The equilibrium condition defining a linear stress says that the sum of normals $\mathbf{n}(F, C)$ weighted by $s(C)$ lies in the linear span of $F$. Higher-dimensional stresses were also considered by Tay, White, and Whiteley [23] and Rybnikov [18]. Our terminology is in good agreement with terminology in these papers.

If $f_k$ denotes the number of simplexes of dimension $k$ in $K$, then "magic" numbers $g_k$ and $h_k$ are defined as follows

$$g_k(K, d) = \sum_{j=-1}^{k-1} (-1)^{k+j-1} \binom{d-j}{d-k+1} f_j$$

$$h_k(K, d) = \sum_{j=0}^{k} (-1)^{j+k} \binom{d-j}{d-k} f_{j-1}$$

For a generic realization in $\mathbb{R}^{d+1}$ of a simplicial homology $d$-manifold $\Delta$ with homology groups of a standard sphere the dimension of the space of $k$-stresses is $g_k(\Delta, d+1)$ if $k \leq \lfloor \frac{d+1}{2} \rfloor$, and 0 if $k > \lfloor \frac{d+1}{2} \rfloor$ (see [13]). For a generic realization in $\mathbb{R}^d$ of a simplicial homology $d$-manifold $\Delta$ with homology groups of a standard sphere the dimension of the space of $k$-stresses is $h_k(\Delta, d+1)$ [23, 13].

There is no similar algebraic theory of stresses for non-simplicial manifolds. The main barrier is the absence of an analog of the notion of face ring for non-simplicial complexes.

## 3 Orientability and generalized volumes.

Let $\mathbb{R}^d$ be a Euclidean affine space with a fixed coordinate system. Consider an oriented, simplicial $(d-1)$-manifold $\Delta$ realized in $\mathbb{R}^d$. We introduce a generalized volume function, $Vol_d$, which assumes positive, negative or zero values on such manifolds. In the case where the manifold $\Delta$ bounds a $d$-dimensional body, and the orientation of $\Delta$ is chosen appropriately, $Vol_d(\Delta)$ is the standard Euclidean volume of the body. Let $F = (v_1, \ldots, v_d)$ be an oriented simplex in $\mathbb{R}^d$. We denote by $[\mathbf{v}_1(F) - \mathbf{p}, \ldots, \mathbf{v}_d(F) - \mathbf{p}]$ the matrix whose rows are $d$-vectors pointing from point $\mathbf{p} \in \mathbb{R}^d$ to the vertices of $F$.

**Definition 3.1** *Let $\Delta$ be a closed oriented simplicial manifold of co-dimension 1 in $\mathbb{R}^d$. Then*

$$Vol_d(\Delta) = \frac{1}{d!} \sum_{F \subset \Delta} det[\mathbf{v}_1(F) - \mathbf{p}, \ldots, \mathbf{v}_d(F) - \mathbf{p}]$$

*where the summation ranges over all oriented $(d-1)$-faces of $\Delta$.*

It is well known that the generalized volume does not depend on the choice of point $\mathbf{p}$. That is why it is normally written for $\mathbf{p} = \mathbf{0}$. The above formula can be rewritten as

(1) $$Vol_d(\Delta) = \frac{1}{d} \sum_{F \in \Delta} d(\mathbf{p}, aff(F)) Vol_{d-1}(F, \mathbf{p})$$

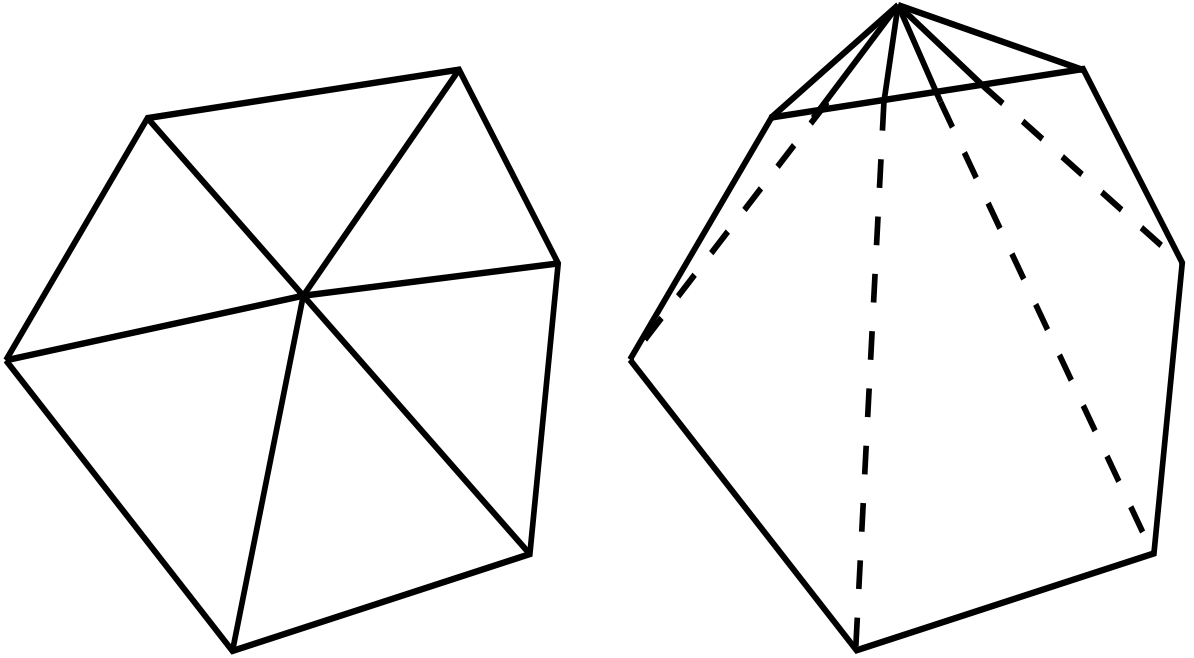

Figure 4: Two realizations of a star

where $d(\mathbf{p}, aff(F))$ stands for the distance between $\mathbf{p}$ and the hyperplane spanned by face $F$. The generalized $(d-1)$-volume $Vol_{d-1}(F, \mathbf{p})$ is computed with respect to the orientation of $aff(F)$ induced by vector $\mathbf{v}_i(F) - \mathbf{p}$ ($i$ is arbitrary), i.e. with respect to an orthonormal coordinate frame $[\mathbf{e}_1 \ldots \mathbf{e}_{d-1}]$ in $aff(F)$ such that $[\mathbf{v}_i(F) - \mathbf{p}, \mathbf{e}_1 \ldots \mathbf{e}_{d-1}]$ is positively oriented in $\mathbb{R}^d$.

Let $\mathbb{S}^{d-1}$ be an oriented simplicial sphere, and let $D$ be a cell-decomposition of $\mathbb{S}^{d-1}$ which is the result of an amalgamation of some of the simplexes of into blocks $\mathbb{S}^{d-1}$ (see Section 1 and [16, 21]). Consider a realization of the simplicial complex $\mathbb{S}^{d-1}$ in $\mathbb{R}^d$ such that each block lies in the affine span of its vertex set. For example, a block can be realized as a convex polytope partitioned into simplexes or as a simplicial star with self-intersections (see Figure 4). Then $Vol_d(\mathbb{S}^{d-1})$ does not depend on the positions of the baricenters of the blocks of all dimensions greater than 0. It can be shown by induction in $d$. The case of $d = 1$ is obvious. The induction step follows from an application of Formula 1.

In Section 6 we will use the following observation.

**Remark 3.2** *Let $B$ be a $d$-dimensional cell-complex such that the closures of all its faces, including $B$, are cones over homology spheres. In other words $B$ is a homology ball. An example of such complex would be a convex polytope. Suppose a baricentric triangulation of $B$ is realized in $\mathbb{R}^d$ so that the affine dimension of the vertex set of each cell of $B$ equals to the dimension of this cell (the cell structure of a convex polytope would serve as a simple example). Then the generalized volume of this simplicial sphere does not depend* *on the positions of the baricenters of its faces provided the baricenter of each face of $B$
lies in the affine span of the vertex set of this face. We call the realizations of these
baricenters in $\mathbb{R}^d$ virtual baricenters.*

For discussion of the algebraic properties of the generalized volume $Vol_d(\mathbb{S}^{d-1})$ as function of the edge lengths see [5].

## 4 Combinatorial dual graph and reciprocals

Let $F(V, E)$ be a framework realized in $\mathbb{R}^2$, and assume that graph $(V, E)$ can be regarded as the 1-skeleton of a spherical complex $\Delta$. Suppose that this framework is in a state of static equilibrium. Consider a vertex of $(V, E)$. The sum of vectors of stresses applied to this vertex is equal to zero. Therefore, when rotated on 90° clockwise they form a polygon (self-intersecting in general). It was first noticed by Maxwell (and proved by Whiteley [26]) that the positions of rotated edges of $F(V, E)$ can be adjusted so that they form a reciprocal graph (often called simply *reciprocal*). Each edge of this reciprocal corresponds to an edge of $F(V, E)$ and each vertex to a cell of $\Delta$. One can introduce a similar notion for piecewise-linear realizations of $d$-manifolds in $\mathbb{R}^d$ (for more information see [2, 18, 25, 10]). In this section we will explore connections between the $d$-stresses and the generalization of Maxwell reciprocal for $d$-manifolds.

The *combinatorial dual graph* $\mathcal{G}(\mathcal{M}^d)$ of a manifold $\mathcal{M}^d$ is defined as follows. The vertices of $\mathcal{G}$ are $d$-cells of $\mathcal{M}^d$, and the edges of $\mathcal{G}$ are internal $(d-1)$-cells of $\mathcal{M}^d$. Two vertices share an edge if and only if the corresponding $d$-cells are adjacent.

A *reciprocal* of a piecewise-linear realization $M$ of a manifold $\Delta$ in $\mathbb{R}^d$ is a rectilinear realization $R$ in $\mathbb{R}^d$ of the combinatorial dual graph $\mathcal{G}(\Delta)$ such that the edges of $R$ are perpendicular to the corresponding facets. If none of the edges of a reciprocal collapses into a point, the reciprocal is called non-degenerate. Reciprocals were originally considered by Maxwell [14] in connection with stresses in plane frameworks. He and almost at the same time L. Cremona [11] noticed that convex reciprocals corresponded to convex liftings of planar cell-complexes. Reciprocals were later studied in [2, 8, 9, 26, 19]. Crapo and Whiteley gave an explicit treatment of the theory of reciprocals, stresses and liftings for 2-manifolds in [8, 9, 10].

To illustrate the concept of reciprocal let us consider the case where the realization $M$ is an embedding. Let $v(C_1)$ and $v(C_2)$ be vertices of a reciprocal $R$ corresponding to adjacent $d$-cells $C_1$ and $C_2$. Call the edge $[v(C_2)v(C_1)]$ *properly oriented* if $\mathbf{v}(C_2) - \mathbf{v}(C_1)$ is cooriented with an outer normal to $C_1$ at the facet shared with $C_2$. Otherwise call $[v(C_2)v(C_1)]$ improperly oriented. A hexagonal reciprocal for the embedded star of a vertex in a 2-manifold is shown on Figure 5. One can see that edges $ef$, $cd$ are improperly oriented, and edges $ab$, $cb$, $de$, and $fa$ are properly oriented). If all edges of $R$ are properly oriented $R$ is called a convex reciprocal (since the cycles of $R$ corresponding to the stars of the $(d-2)$-cells are convex in this case).We refer to reciprocals of stars of the manifold as local reciprocals.

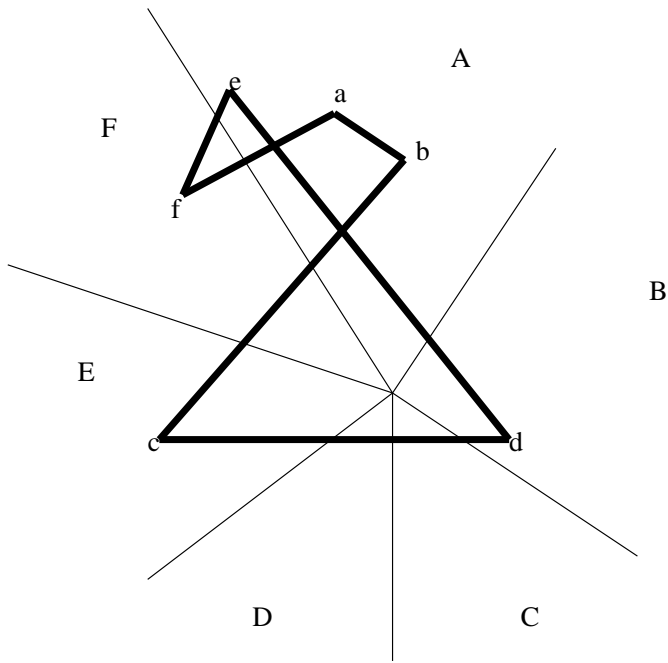

Figure 5: Non-convex reciprocal

Evidently, reciprocals with one fixed vertex form a linear space. Denote it by $Rec(M)$. If $M$ is an embedding, then convex reciprocals form a cone $CRec(M)$ in the linear space $Rec(M)$. The following theorem by Rybnikov [18] explains connections between reciprocals and stresses in the case of general dimension. We will utilize this theorem in the proof of our main theorem from Section 6.

**Theorem 4.1** *Let $M$ be a PL-realization of a homology $d$-manifold $\Delta$ in $\mathbb{R}^d$ with trivial first homology group over $\mathbb{Z}/2\mathbb{Z}$. Then there is an isomorphism between $Stress_d(M)$ and $Rec(M)$. Non-zero coefficients of stresses correspond to non-vanishing edges of a reciprocal. If $M$ is an embedding of $\Delta$ into $\mathbb{R}^d$, then one can interpret properly oriented edges as corresponding to tensed facets, and improperly oriented edges as corresponding to compressed facets.*

Let $B$ be a $d$-dimensional cell-complex which is the cone over a homology sphere (not necessarily simplicial). Obviously, $B \setminus \partial B$ can be regarded as a star $St$. Let $R$ be a reciprocal for $St$ and denote by $R(C)$ a sub-reciprocal of $R$ corresponding to a face $C \in St$. The vertex set of $R$ is a realization of the vertex set of a complex dual to $St$. Denote it by $St^*$. For each cell $C$ ($k = 1 \leq dim(C) \leq d$) of $St^*$ choose an arbitrary point $vbc(C, R)$ on each plane $aff(R(C))$, and call it the *virtual baricenter* of $R(C)$. The vertices of $R$ and the points $vbc(C, R)$, $k = 1 \leq dim(C) \leq d$ define a piecewise-linear realization of $St^*$. We know from Remark 3.2 that if a baricentric triangulation of $St^*$ is realized in $\mathbb{R}^d$ so that the affine dimension of the vertex set of each cell of $St^*$ equals to the dimension of this cell (the cell structure of a convex polytope would serve as a simple example), then the generalized volume of oriented simplicial sphere $\partial\, St^*$ does not depend on the positions of the virtual baricenters of its faces provided the virtual baricenter of each face of $St^*$ lies in the affine span of the vertex set of this face. We can sum up this observation in the following proposition which will be of great use in the following section.

**Proposition 4.2** *Let $R$ be a reciprocal for an oriented $d$-dimensional star $St$ realized in $\mathbb{R}^d$. Then the generalized volume $Vol_d(R)$ is well defined.*

## 5  Minkowski theorem and stresses

In this section we give an application of a well-known Minkowski theorem to stresses on polyhedral partitions of $\mathbb{R}^d$.

**Theorem 5.1** *(Minkowski) Let $P$ be a convex polytope in $\mathbb{R}^d$, and denote by $\{\mathbf{n}(F)\}$ the inner unit normals to facets of $P$. Then*

$$\sum_{F \subset \partial P} vol_{d-1}(F)\, \mathbf{n}(F) = \mathbf{0}.$$

If we choose a (combinatorial) orientation for $P$ and denote by $Vol_{d-1}(F, \mathbf{n}(F))$ the generalized volume of an oriented facet $F$ with respect to the orientation of $aff(F)$ induced by $\mathbf{n}(F)$, then the above formula can be rewritten as

$$\sum_{F \subset \partial P} Vol_{d-1}(F, \mathbf{n}(F))\, \mathbf{n}(F) = \mathbf{0}.$$

Notice that in the last formula the directions of normals $\mathbf{n}(F)$ need not agree, since $Vol_{d-1}(F)$ is computed with respect to the orientation induced by $\mathbf{n}(F)$. Flipping the normal changes the sign of $Vol_{d-1}(F)$.

Let $St(v)$ be the star of a vertex $v$ of a polyhedral partition of $\mathbb{R}^d$.

**Definition 5.1** *A dual convex polytope for $St(v)$ is a $d$-dimensional polytope $D(St(v))$ in $\mathbb{R}^d$ satisfying the following conditions.*

*1) There is a one-to-one correspondence $\mathcal{I}$ between the $m$-dimensional faces of $D(St(v))$ and the $(d-m)$-dimensional faces of $St(v)$ ($0 \leq m \leq d$).*

*2) If $D^s \subseteq D^t$ are faces of $D(St(v))$ corresponding to faces $F^{d-s}$ and $F^{d-t}$ of $St(v)$, then $F^{d-t} \subseteq F^{d-s}$. In other words the mapping $\mathcal{I}$ induces an isomorphism between the face lattices of $D(St(v))$ and $St(v)$.*

*3) For $0 \leq m \leq d$ each $m$-dimensional face of $D(St(v))$ is perpendicular to the corresponding $(d-m)$-dimensional face of $St(v)$.*

*4) $Sk^1(D(St(v)))$ is a convex reciprocal graph for the star $St(v)$ (see Section 4).*

The convexity of the dual polytope immediately follows from Conditions 1 - 4. Suppose that there is a $d$-tension on St(v) (all coefficients of $d$-stresses are strictly positive). By results of [18] and [19] there is a convex polytope $D$ dual to $St(v)$. $D$ is uniquely determined by this tension up to translation. By the Minkowski theorem cited above the sum of facet normals of a convex polytope scaled by the facet volumes is zero. Therefore, one can interpret the volumes of $m$-faces, $1 \leq m \leq d-1$ of $D$ as coefficients of $(d-m+1)$-stresses on $(d-m)$-dimensional cells of $St(v)$. Thus a $d$-tension on the star $St(v)$ induces an $(d-m)$-tension on $St(v)$, $1 \leq m \leq d-1$. It is easy to see that the constructed mappings are polynomial.

**Proposition 5.2** *Let $\Delta$ be a cell-decomposition of a polyhedral region in $\mathbb{R}^d$. For $k = 1, \ldots, d-1$ there is a polynomial mapping of degree $d-k+1$ from the cone of $d$-tensions of $\Delta$ to the cone of $k$-tensions of $\Delta$. An all non-zero tension is always mapped to an all non-zero tension.*

By construction, in the case of embedding a $d$-tension is mapped to a 2-tension on the 1-skeleton of the manifold.

**Corollary 5.3** *Let $G$ be the 1-skeleton of a cell-decomposition $\Delta$ of $\mathbb{R}^d$ by convex polyhedra. If there is a convex surface which projects onto $\Delta$, then $G$ supports a positive equilibrium stress at all edges, and therefore is an infinite spider web.*

It turns out that the mappings from Proposition 5.2 can be extended from the cone of tensions to all the space of $d$-stresses, and the above construction can be carried out for arbitrary piecewise-linear realizations of orientable $d$-manifolds (not necessary embeddings). In order to formally establish this, we will need the concept of generalized volume introduced in Section 3.

The Minkowski theorem can be formulated for simplicial spheres arbitrarily realized in $\mathbb{R}^d$, and as we will see in Section 6 even for a large class of non-simplicial spheres realized in $\mathbb{R}^d$ with self-intersections. Let $\Delta$ be an oriented simplicial manifold realized in $\mathbb{R}^d$. For each oriented $(d-1)$-simplex $F$ pick a unit normal vector $\mathbf{n}(aff(F))$, and let $Vol_{d-1}(F, \mathbf{n}(aff(F)))$ be the generalized volume of $F$ computed in $aff(F)$ equipped with an orientation induced by $\mathbf{n}(F)$. We need the following lemma.

**Lemma 5.4**

$$\sum_{F \subset \Delta} Vol_{d-1}(F, \mathbf{n}(aff(F))) \, \mathbf{n}(aff(F)) = \mathbf{0}.$$

*Proof.* The orientation of $\Delta$ induces an orientation on a cone with $\Delta$ as base. Thus if $F_1$ and $F_2$ are two adjacent $(d-1)$-faces of $\Delta$, the orientations of the cone over their common facet are opposite. Therefore the above formula can be rewritten as

$$\sum_{F \subset \Delta} \sum_{s^{d-1} \subset \partial_0 \cdot F} Vol_{d-1}(s^{d-1}, \mathbf{n}(aff(s^{d-1}))) \, \mathbf{n}(aff(s^{d-1})) = \mathbf{0}.$$

where $s^{d-1}$ stands for a facet of the cone $\mathbf{0} \cdot F$, and $\mathbf{n}(aff(s^{d-1}))$ is an arbitrary unit normal to hyperplane $aff(s^{d-1})$. Applying Minkowski theorem to each $d$-simplex $\mathbf{0} \cdot F$ we get the required formula. □

The interplay between stresses and volumes in the case of convex polytope was also described by McMullen [17] and Lee [13].

## 6  A natural trace of $d$-stresses in lower dimensions.

**Remark 6.1** *Let $\Delta$ be an orientable homology $(d-1)$-manifold in Euclidean space of dimension $d$. An orientation of $\Delta$ induces the orientation of normals to $\Delta$ at the cells of maximal dimension by the following rule. Let $(v_1(S), \ldots, v_d(S))$ be an oriented simplex of $\Delta$. If frame $[\mathbf{v}_1(S), \ldots, \mathbf{v}_d(S)]$ is positively oriented, then the corresponding normal to $\Delta$ at $S$ has positive scalar product with all these vectors. Conversly, a consistent choice of the field of normals to $\Delta$ at their simplexes of maximal dimension determines orientation of $\Delta$ (e.g. outer normals for a convex polytope; see Figure 6).*

In the case of an orientable $d$-manifold it is possible to fix the orientation of cells so that they form a $d$-cycle. By the above remark such orientation of cells induces the orientation of frames of normals corresponding to flags of cells. Thus, if $\Delta$ is an orientable $d$-manifold in $\mathbb{R}^d$ and the orientations of $d$-cells are picked up in such a way that it turns $\Delta$ into a $d$-cycle, any two flags of equal length having $d$-cells as maximal elements and distinct at only one position have corresponding frames of opposite orientation. A face-to-face partition of $\mathbb{R}^d$ provides a transparent example. Each of the two possible orientations of the partition correspond to either flags of inner normals or to the flags of outer normals.

**Theorem 6.2** *Let $\Delta$ be an orientable homology $d$-manifold realized in $\mathbb{R}^d$. Then for $k = 1, \ldots, d-1$ there is a polynomial mapping $\mathfrak{p}_k$ of degree $d-k+1$ from the $Stress_d(\Delta)$ to $Stress_k(\Delta)$.*

*Proof.* For a cell-decomposition of a homology manifold there is so-called dual cell-decomposition (also called dual block decomposition). Consider the baricentric triangulation $T(\Delta)$ of the original cell-decomposition. Each cell of the original decomposition is a simplicial star in the baricentric triangulation. All $(d-k)$-simplexes of this triangulation which share the baricenter of a $k$-cell $c$ form the dual cell (also called block) for $c$. This dual cell is a homology $(d-k)$-disk. The boundary of the dual cell is a homology sphere (for more details on the geometrical duality in homology manifolds see [16, 21]).

Let $v$ be a vertex of $\Delta$, and let $D_v$ be the $d$-dimensional cell (block) corresponding to $v$ in the dual decomposition of $\Delta$. Obviously, the boundary of $D_v$ is the link $Lk(v)$ of $v$. Each $k$-simplex of the baricentric triangulation of $D_v$, $k = 0, \ldots, d-1$, can be regarded as the result of the $(k-1)$-fold iterative coning starting from vertex $v$.

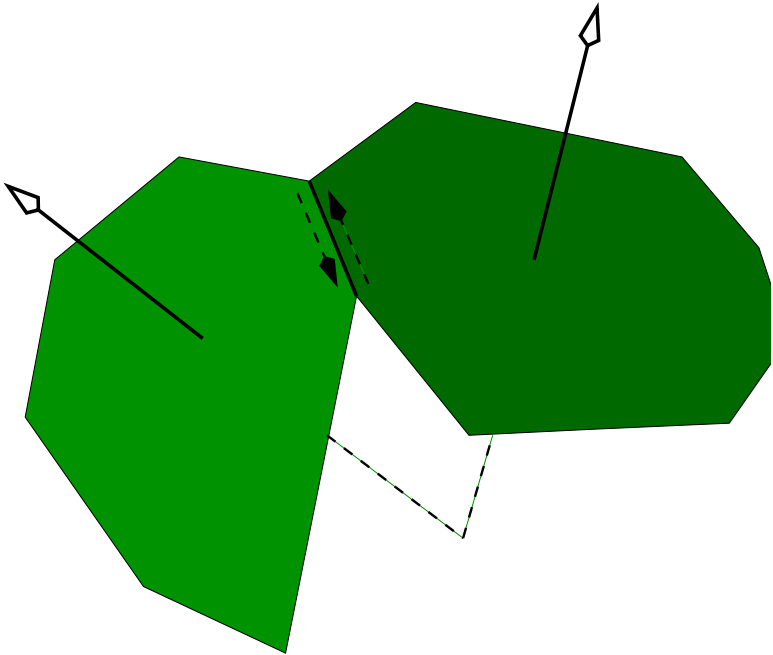

Figure 6: Orientation

Each cell of $\Delta$ or $\Delta^*$ is itself an orientable homology manifold, namely a homology disk. Thus, an orientation of triangulation $T(\Delta)$ induces in a natural way orientation on cell-complexes $\Delta$ and $\Delta^*$.

Let $R$ be a (Euclidean) reciprocal for $St(v) \subset \Delta$ (see Section 4). By Theorem 4.1, the linear space of $d$-stresses on $St(v)$ is naturally isomorphic to the space of reciprocals with one fixed vertex. It turns out that one can introduce the notion of generalized "$k$-volume" ($k = 0, \ldots, d$) for the sub-reciprocals of $R$, corresponding to the stars of cells of $St(v)$ (we refer to them as "faces" of R). It is natural to call this function $k$-volume, because when a reciprocal can be regarded as the vertex set of a convex $k$-polytope, the absolute value of this function is equal to the $k$-volume of the polytope. We keep the same notation for the $k$-volume of a reciprocal that we used for generalized volumes, i.e. $Vol_k$.

Let $C^{d-k}$ be a $d - k$-cell from the (open) star of $v$. Obviously $St(C^{d-k}) \subset St(v)$. The subset $R(C^{d-k})$ of $R$ corresponding to this star spans an affine $k$-plane perpendicular to $C^{d-k}$. $R(C^{d-k})$ can be regarded as a realization of the vertex set of a cell of $\Delta^*$ dual to $C^{d-k}$. Thus it makes sense to talk about the (combinatorial) orientation of $R(C^{d-k})$. Recall that a $k$-cell of the dual decomposition corresponds to a $(d-k)$-cell of $\Delta$. Choose a flag of full length in $C^{d-k}$. It corresponds to some simplex $s$ of $T(\Delta)$ whose vertex set is the "baricenters" of the flag cells. The iterative coning of $C^k$ with vertices of $s$ is a cell from an amalgamation of triangulation of the star of $v$ in $\Delta^*$ into (non-simplicial, in general) blocks of form $v_0 \cdots (v_{d-k} \cdot Dual(C^k))$ constructed by successful coning of $C^k$. An orientation of $T$ induces an orientation on $v_0 \cdots (v_{d-k} \cdot C^k)$. Therefore the choice of flag determines an orientation for $C^k$.

A flag of faces of $C^{d-k}$ corresponds to an ordered $(d-k)$-tuple of normals to the faces of $C^{d-k}$. Denote it by $[N]$. This $(d-k)$-tuple induces an orientation of affine subspace spanned by $R(C^{d-k})$ by the following rule. A frame $F$ in $aff(R(C^{d-k}))$ is said to be cooriented with the frame $[N]$ if $[N, F]$ is cooriented with the coordinate frame of $\mathbb{R}^d$. Therefore $Vol_{d-k}(R(C^{d-k})$ is well defined provided a flag of cells in $C^{d-k}$ (see Section 3) has been fixed. We have to show that $Vol_{d-k}(R(C^{d-k}))$ does not depend on the choice of flag in $C^{d-k}$. It is enough to show that for two flags in $C^{d-k}$ that differ in one position the $Vol_{d-k}(R(C^{d-k}))$ is the same, since any two flags in $C^{d-k}$ can be connected by a sequence of alterations. Obviously, two flags that differ in one position induce opposite combinatorial orientations on $R(C^{d-k})$. But on the other hand it means that the $(d-k)$-tuples of vectors corresponding to these flags have opposite orientations. Thus the generalized $k$-volume of $R(C^{d-k})$ is well defined and does not depend on the choice of a flag of faces in $C^{d-k}$.

Let $s$ be a $d$-stress on $\Delta$. Since the star of a $(d-k)$-cell of $\Delta$ is a homology $d$-disk, a $d$-stress restricted to the star of a vertex generates a $k$-dimensional reciprocal for this star (see Section 4). The distance between two vertices of the reciprocal corresponding to two adjacent $d$-cells equals the absolute value of stress on their common facet. Let $R(C)$ be the reciprocal of the star of a $(d-k)$-cell $C$ corresponding to the stress $s$. Let us interpret $Vol_k(R(C))$ as the value of $(d-k+1)$-stress on $C$ (recall that $(d-k)$-cells bear $(d-k+1)$-stresses). We have to check the equilibrium condition at every $(d-k-1)$-cell of $\Delta$. Let $F$ be a $(d-k-1)$-cell of $\Delta$. Construct the reciprocal $R(F)$ for $St(F)$ corresponding to the $d$-stress $s$. Notice that if $F \subset C$, then the sub-reciprocal or $R(F)$ corresponding to the star of $C$ coincides with $R(C)$ (up to translation). Let $\mathbf{n}(F, C)$ denotes the fixed unit normal to $C$ at $F$ whose orientation is induced by the orientation of $\Delta$ as it was explained in the beginning of this section. In the case where $\Delta$ is embedded into $\mathbb{R}^d$ we can think of $\mathbf{n}(F, C)$ as of inward unit normal.

$$\sum_{\{C|\ F \subset C\}} Vol_k\left(R(C), \mathbf{n}(F, C)\right) \mathbf{n}(F, C) = \sum_{\{R(C)|\ F \subset C\}} \sum_{S \subset C} Vol_k(S, \mathbf{n}(F, C)) \mathbf{n}(F, C)$$

where $S$ is an oriented $(d-k)$-simplex from a baricentric triangulation of $R(C)$ arbitrarily realized in $aff(R(C))$. By Minkowsi theorem the last quantity is always zero. □

One should notice that the orientability of $\Delta$ is essential for our construction. Only in the case of orientable manifold the edges of a reciprocal can be separated into properly oriented and improperly oriented.

Since the generalized $(d-k+1)$-volume of $R$ is a homogeneous polynomial of degree $d-k+1$ in the (oriented) lengths of the edges of $R$, and the absolute values of the edges of $R$ equal to the absolute values of corresponding $d$-stresses (see Section 3), the constructed mappings $\mathfrak{p}_k$ from $Stress_d(\Delta)$ to $Stress_k(\Delta)$, $k = 1, \ldots, d-1$ are polynomial of degree $d-k+1$. The coefficients of these polynomials depend on geometry of $\Delta$. Mapping $\mathfrak{p}_k$ can also be regarded as the restrictions of certain rational $\mathbb{R}^{f_{k-1}}$-valued function $\mathfrak{m}_k$ defined on Euclidean space of dimension $f_{d-k+1}$ to the linear subspace of $d$-stresses. The values of $k$-stresses are taken from the intersection of the image of $\mathfrak{m}_k$ and the subspace $Stress_k(\Delta)$ of $\mathbb{R}^{f_{k-1}}$. According to Connelly, Sabitov, and Waltz [5, 31] the volume of an orientable simplicial $d$-manifold for $d = 3, 4$ is a function of the edge lengths only. It means that when $d = 3, 4$ and $k = d - 2, d - 3$, mapping $\mathfrak{p}_k$ considered as a multivaried function on $\mathbb{R}^N$, $N = f_{d-k+1}$, for each coordinate in the image is a ratio of polynomials with integer coefficients that do not depend on the realization of the manifold $\Delta$ in $\mathbb{R}^d$. It would be interesting to know if there are any implications of this fact for the algebraic geometry of our mappings $\mathfrak{p}_k$.

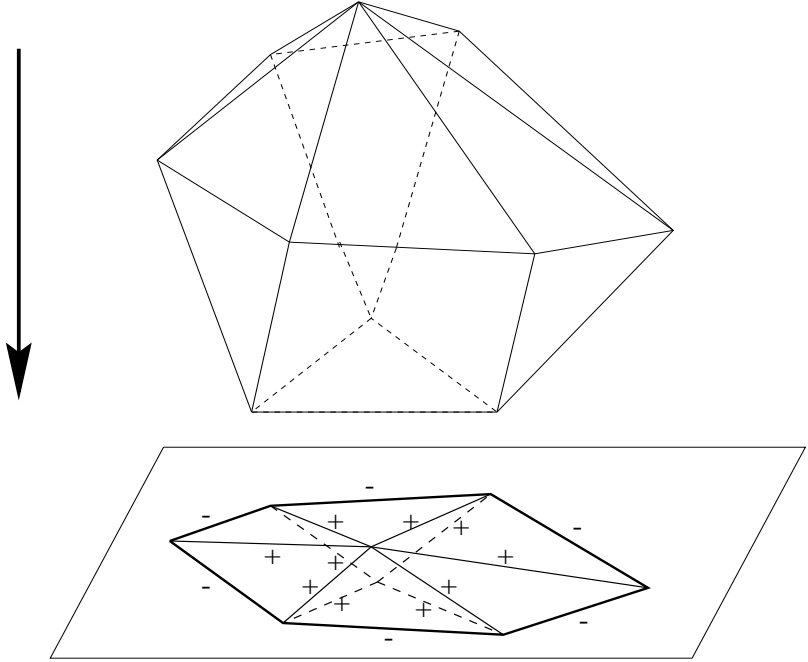

Figure 7: Maxwell convex stress

By construction, in the case of an embedding a $d$-tension is mapped to a 2-tension on the 1-skeleton of the manifold.

**Corollary 6.3** *Let $G$ be the 1-skeleton of a decomposition $\Delta$ of $\mathbb{R}^d$ by convex polyhedra. If there is a convex surface which projects onto $\Delta$, then $G$ is a spider web.*

Maxwell [14, 15] discovered the "convex self-stress" induced by projection of a convex polytope on the plane (see Figure 7).

**Theorem 6.4** *The vertical projection of a strictly convex polyhedron, with no faces vertical, produces a plane framework with a self-stress that is negative on the boundary edges and positive on all edges interior to this boundary polygon.*

Now we can formulate a partial analog of Maxwell theorem on convex self-stresses and projections of spatial polyhedra. It immediately follows from our main theorem.

**Theorem 6.5** *Let $P^4$ be a strictly convex polytope in $\mathbb{R}^4$ without vertical faces, and let $G$ be the projection of $Sk^1(P^4)$ onto $\mathbb{R}^3 \subset \mathbb{R}^4$. Then $G$ supports a self-stress $s$ which is positive on all edges of $G$ that belong to the interior of the projection. If all the edges of $P^4$ that project on the boundary of the projection are incident to exactly three 3-cells of $P^4$, then in addition $s$ is negative on all edges of $G$ that belong to the boundary of the projection.*

*Proof.* Using our main theorem, let us construct the mapping $\mathfrak{p}_2 : Stress_3 \to Stress_2$ for the realization of our polytope $P^4$ in $\mathbb{R}^3$ induced by the vertical projection. Obviously, since the upper and the lower lids are convex, the reciprocals for the "interior" edges are convex (1-skeletons of convex polytopes) and have volumes of the same sign. The reciprocals of the "boundary" edges need not be convex; however if a boundary edge has a simplicial reciprocal, its volume ought to have the sign opposite to signs of the volumes of the reciprocals of the interior edges. □

Recall than Maxwell correspondence states also that any equilibrium stress can be interpreted as one induced by the projection of a spatial polytope. On the CMS winter meeting of 1998 R. Connelly and W. Whiteley asked if the following conjecture is true for our correspondence.

**Conjecture 6.6** *Let $M^3$ be a homology sphere realized in $\mathbb{R}^4$ and let $s_2$ be a self-stress (2-stress) on the 1-skeleton of $M^3$. There is a 3-stress $s_3$ on $M^3$ such that $\mathfrak{p}_2(s_3) = s_2$.*

As it was mentioned in the introduction, the generic realization of the boundary of the 4-dimensional cross-polytope $O_4$ provides a counterexample. According to Lee [13] $dim(Stress_3(O_4)) = 4$, but $dim(Stress_2(O_4)) = 6$ (self-stresses on a framework are 2-stresses). Since the mappings are algebraic the image of the space of 3-stresses cannot cover the space of 2-stresses. It would be interesting to give a geometric interpretation of those 2-stresses that can be interpreted as images of 3-stresses under the above mapping.

A cell-decomposition of a closed $d$-manifold is called $k$-primitive if the star of each $k$-cell has $d - k + 1$ $d$-cells (some authors call 0-primitive decompositions *simple*; our terminology goes back to Voronoi [24]). The meaning of this definition is that in a decomposition of $\mathbb{R}^d$ by convex polyhedra, $d - k + 1$ is the minimal possible number of $d$-cells making contact in a $k$-cell. When a $k$-primitive cell-decomposition of $\mathcal{M}^d$ is assumed to be fixed, we will refer to this $k$-primitive decomposition of $\mathcal{M}^d$ as $k$-*primitive manifold* $\mathcal{M}^d$. If a PL-realization of a sphere $\mathbb{S}^d$ in $\mathbb{R}^d$ can be lifted to a convex polytope in $\mathbb{R}^{d+1}$, then 0-primitive vertices of $\mathbb{S}^d$ correspond to simple vertices of this convex polytope. The notion of $k$-primitive decomposition naturally arises in studies of space-fillers, lattice polytopes and stereohedra. For example, the affine equivalence between space-fillers and Dirichlet domains of lattices was proved by Voronoi only for 0-primitive (simple) tilings. The existence of a lattice Dirichlet domain which is affinely isomorphic to a space-filler $\Pi$ is equivalent to the existence of a $d$-stress with some special symmetries on the lattice tiling $T(\Pi)$ by $\Pi$ (Voronoi)[18, 19, 24] . Since any $(d-3)$-primitive decomposition of $\mathbb{R}^d$ is the projection of a convex surface [18, 19], we have the following corollaries.

**Corollary 6.7** *The 1-skeleton of a $(d-3)$-primitive decomposition of $\mathbb{R}^d$ by convex polyhedra is always a spider web.*

A cell-decomposition of a $d$-manifold is referred to as $k$-*primitive* if the star of each internal $k$-dimensional cell has $d-k+1$ $d$-cells (some authors call 0-primitive decompositions *simple*; our terminology goes back to Voronoi [24]). For decompositions of $\mathbb{R}^d$ by convex polyhedra $d-k+1$ is the minimum possible number of tiles in the star of a $k$-face.

**Corollary 6.8** *Let $M$ be a realization in $\mathbb{R}^d$ of a $(d-3)$-primitive manifold $\Delta$ with trivial $H_1(\Delta, \mathbb{Z}_2)$. Suppose the body $|M|$ of this realization is convex and $M$ is a double cover of int $|M|$. Then the 1-skeleton of $M$ admits a convex self-stress.*

**Conjecture 6.9** *Let $\Delta$ be a simplicial homology $d$-manifold with $H_1(\Delta, \mathbb{Z}_2) = 0$. Then for any generic realization of $\Delta$ in $\mathbb{R}^d$ mappings $\mathfrak{p}_k$, $k = 1, \ldots, d$ have Jacobian of maximal possible rank at almost all points $s \in Stress_d(\Delta)$.*

One can ask about generic properties of the mappings $\mathfrak{p}_k$, only when any realization of $\Delta$ admits small perturbations not changing its combinatorial structure. For instance, this is the case when $d = 2$, when $\Delta$ is simplicial or when $\Delta$ admits a sharp lifting (for details on sufficient conditions for the existence of a sharp lifting see [18]). It is plausible that in these cases, for generic realizations the constructed mappings also have Jacobian of maximal possible rank. It is possible that the above conjecture holds for arbitrary orientable manifolds.

A necessary condition for our theorem is that $dim(Stress_d) \leq dim(Stress_k)$, $k > d$. Below we give a count that demonstrates that this condition holds for $k = 2$ (mapping of a $d$-stress to a self-stresses on the 1-skeleton). The dimension of the space of $d$-stresses on a simplicial $d$-pseudomanifold in $\mathbb{R}^d$ is at least $f_0 - d - 1$ [4] and is equal to $f_0 - d - 1$ if $\Delta$ is a manifold with $H_1(\Delta, \mathbb{Z}_2) = 0$ [18]. By the result of Fogelsanger [12] the 1-skeleton of a generic realization of a $d$-pseudomanifold in $\mathbb{R}^{d+1}$ is statically rigid. It means that $Sk^1(\Delta)$ can resolve any external load in $\mathbb{R}^{d+1}$ (see Introduction). Thus $dim\ Stress_2(\Delta, d+1) = f_1 - (d+1)f_0 + \binom{d+2}{2} = g_2(\Delta, d+1) \geq 0$ (the lower bound theorem for general simplicial pseudomanifolds).

For Conjecture 6.9 to be true, it is necessary that

$$dim\ Stress_2(\Delta, d) \geq dim\ Stress_d(\Delta, d) = f_0 - d - 1.$$

By the lower bound theorem $dim\ Stress_2(\Delta, d) - (f_0 - d - 1) = f_1 - (d+1)f_0 + \binom{d+2}{2} = dim\ Stress_2(\Delta, d+1) = g_2(\Delta, d+1) \geq 0$.